\documentclass[12pt]{article}
\usepackage{amsfonts}
\usepackage{amssymb}
\usepackage{mathtools}
\usepackage{verbatim}
\usepackage{amsmath}
\usepackage{amscd}
\newtheorem{theorem}{Theorem}
\newtheorem{theorem*}{Theorem}
\newtheorem{corollary}{Corollary}
\newtheorem{proposition}{Proposition}
\newtheorem{definition}{Definition}
\newtheorem{definition*}{Definition}
\newtheorem{remark}{Remark}

\begin{document}

\begin{center}
{{\Large \bf On  packing of  Minkowski balls. II
} \\
}
\end{center}

\begin{center}
{\bf Nikolaj M. Glazunov } \end{center}

\begin{center}
{\rm Glushkov Institute of Cybernetics NASU, Kiev, } \\
{\rm Institute of Mathematics and Informatics Bulgarian Academy of Sciences }\\
{\rm  Email:} {\it glanm@yahoo.com }
\end{center} 

\bigskip

{\bf Abstract.}

This is the continuation of the author's ArXiv presentation  `'On  packing of  Minkowski balls. I''.
In section 2 we investigate lattice packings of  Minkowski balls  and   domains. By results of the proof of Minkowski conjecture about the critical determinant we devide the balls and domains on 3 classes: Minkowski, Davis  and  Chebyshev-Cohn.  The optimal lattice packings of the  balls and domains are obtained.  The minimum areas of  hexagons inscribed in the balls and domains and  circumscribed around their are given. Direct limits of direct systems of Minkowski balls and domains and their critical lattices are calculated.   \\



\bigskip

{\bf Thanks.}{ The author is deeply grateful to the Bulgarian Academy of Sciences, the Institute of Mathematics and Informatics of the Bulgarian Academy of Sciences, Professor P. Boyvalenkov for their support.\\
The author was supported by Simons grant 992227.}

\bigskip

\section{Introduction}. 

This is the continuation of the author's ArXiv presentation  `'On  packing of  Minkowski balls. I'' \cite{pmbIv4}.
In section two ``Critical lattices, Minkowski domains and  their optimal packing'' we investigate lattice packings of  Minkowski    domains. By results of the proof of Minkowski conjecture about the critical determinant we devide the domains on 3 classes: Minkowski, Davis  and  Chebyshev-Cohn.  The optimal lattice packings of the   domains are obtained.
Section three gives minimum areas of  inscribed and circumscribed hexagons of  Minkowski,  Davis  and  Chebyshev-Cohn domains.
In section four presents preliminary notions connected with two-dimensional Banach spaces which define Minkowski norms
\[
   |x|^p + |y|^p, \; p \in {\mathbb R}, p \ge 1,
\]
on  real plane ${\mathbb R}^2 = (x,y)$.
 Results of sections two, three and four  lead to  algebro-geometric structures  in the framework of Pontryagin duality theory and its extensions. Direct limits of direct systems of Minkowski balls and domains and their critical lattices are calculated.

\section{Critical lattices, Minkowski domains and  their optimal packing}

   Let $D$ be a fixed bounded symmetric about origin  convex body ({\it centrally symmetric convex body} for short) with volume
$V(D)$.
When considering packing problems for such $D$, it does not matter whether we consider such $D$ with or without a boundary \cite{Cassels,lek}.

\begin{proposition} \cite{Cassels,lek}.
\label{p1}
 If $D$ is symmetric about the origin and convex, then $2D$ is convex and symmetric
 about the origin.
\end{proposition}

 \begin{corollary}
 \label{cor1}
 Let $m$ be integer $m \ge 0$ and $n$ be natural greater $m$. 
If $2^m D$ centrally symmetric convex body then $2^n D$ is again centrally symmetrc convex body.
\end{corollary}
{\bf Proof.} Induction. \\

  We will called $2^m D$ the $2^m$ {\it  dubling} of $D$.

 Definitions of Minkowski, Davis  and  Chebyshev-Cohn balls are given in \cite{pmbIv4}.
 Now we extend their on respective domains.
 
Recall that this division of general Minkowski balls (\ref{eq1}) is based on the results
 of the proof of the Minkowski conjecture \cite{Mi:DA,D:NC,Co:MC,W:MC,GGM:PM,Gl4} and 
the Corollary \ref{cor1}.
  
 For ease of reading, we present the main result of \cite{GGM:PM}.

\begin{theorem}
\cite{GGM:PM}
\label{ggm}
$$\Delta(D_p) = \left\{
                   \begin{array}{lc}
    \Delta(p,1), \; 1 \le p \le 2, \; p \ge p_{0},\\
    \Delta(p,\sigma_p), \;  2 \le p \le p_{0};\\
                     \end{array}
                       \right.
                           $$
here $p_{0} \approx 2.5725$ is a real constant that is defined uniquely by the conditions
$\Delta(p_{0},\sigma_p) = \Delta(p_{0},1)$; note that $2,57 \le p_{0}  \le 2,58. $
\end{theorem}
\begin{remark}
We will call $p_{0}$ the Davis constant.
\end{remark}

\subsection{Domains.}

We consider the following classes of balls (see \cite{pmbIv4}) and domains.

\begin{itemize}
  
\item  {\it  Minkowski domains (MD)}:    $2^m D_p$,  integer $m \ge 1$, for $1 \le p<2$;

\item {\it Davis domains (DD)}: $2^m D_p$,  integer $m \ge 1$, for $p_{0} > p \ge 2$;


\item {\it Chebyshev-Cohn domains (CCD)}: $2^m D_p$,  integer $m \ge 1$,  for $ p \ge p_{0}$;

\end{itemize}

\begin{remark}
\label{r2}
Sometimes, when it comes to using $p$ values that include scopes of different domains, we will use the term Minkowski domains for definitions of different types of domains, specifying
  name when a specific $p$ value or a range of $p$ values containing a single type of domains is specified.
\end{remark}  

Recall from \cite{pmbIv4} two propositions:


\begin{proposition}
\label{p2}
If $\Lambda$ is the critical lattice of the Minkowski ball $D_p$ than the sublattice $\Lambda_2$  of index two of the critical lattice is the critical lattice of $2 D_p$.
\end{proposition}

\begin{proposition} \cite{Cassels,lek}.
\label{p3}
 The dencity of a $(D_p, \Lambda)$-packing is equal to $V(D_p)/d(\Lambda)$ and it is maximal if $\Lambda$ is critical for  $2 D_p$.
\end{proposition}

\begin{proposition}
\label{p7}
Let $m$ be  integer, $m \ge 1$.
If $\Lambda$ is the critical lattice of the 
ball $D_p$ than the sublattice $\Lambda_{2^m}$  of index $2^m$ of the critical lattice is the critical lattice of $2^m D_p$.
\end{proposition}
{\bf Proof.} Since the Minkowski ball $D_p$ is symmetric about the origin and convex, then $2^m D_p$ is convex and symmetric
 about the origin by Corollary \ref{cor1}.
 
When parametrizing admissible lattices $\Lambda$ having three pairs of points on the boundary of the ball $D_p$,
the following parametrization is used \cite{Co:MC,GGM:PM,Gl4}:
 \begin{equation}
\label{eq2}
\Lambda =  \{((1 + \tau^{p})^{-\frac{1}{p}}, \tau(1 + \tau^{p})^{-\frac{1}{p}}), (-(1 + \sigma^p)^{-\frac{1}{p}}, \sigma(1 + \sigma^p)^{-\frac{1}{p}})\}
\end{equation}
 where 
 \[
 0 \le  \tau < \sigma , \; 0 \le \tau \le \tau_p.
  \]
   $\tau_{p}$ is defined by the
equation $ 2(1 - \tau_{p})^{p} = 1 + \tau_{p}^{p}, \; 0 \leq
\tau_{p} < 1. $ 
\[
  1 \le\sigma \le \sigma_p, \; \sigma_p = (2^p - 1)^{\frac{1}{p}}.
\]
Admissible lattices of the form (\ref{eq2}) for $2^m$ {\it  dubling} Minkowski domains $2^m D_p$ have a representation of the form 
 \begin{equation}
\label{eq3}
\Lambda_{2^m D_p} =  \{2^m ((1 + \tau^{p})^{-\frac{1}{p}}, 2^m \tau(1 + \tau^{p})^{-\frac{1}{p}}), (-2^m (1 + \sigma^p)^{-\frac{1}{p}},      2^m \sigma(1 + \sigma^p)^{-\frac{1}{p}})\}
\end{equation}
  Hence the Minkowski-Cohn moduli space for these admissible lattices has the form
 \begin{equation}
 \Delta(p,\sigma)_{2^m D_p} = 2^{2m} (\tau + \sigma)(1 + \tau^{p})^{-\frac{1}{p}}
  (1 + \sigma^p)^{-\frac{1}{p}}, 
 \end{equation}
in the same domain
 $$ {\mathcal M}: \; \infty > p > 1, \; 1 \leq \sigma \leq \sigma_{p} =
 (2^p - 1)^{\frac{1}{p}}, $$
  
Consequently, the critical determinants of $2^m$ {\it  dubling}
  balls $D_p$ have a representation of the form
\begin{equation}
\label{cdde0}
 {\Delta^{(0)}_p}(2^m D_p) = \Delta(p, {\sigma_p})_{2^m D_p} =  2^{m-1}\cdot {\sigma}_{p},   \; {\sigma}_{p} = (2^p - 1)^{1/p},
 \end{equation}
 \begin{equation}
\label{cdde1}
 {\Delta^{(1)}_p}(2^m D_p)  = \Delta(p,1))_{2^m D_p} = 4^{m - \frac{1}{p}}\frac{1 +\tau_p }{1 - \tau_p},   \;  2(1 - \tau_p)^p = 1 + \tau_p^p,  \;  0 \le \tau_p < 1. 
 \end{equation}
  And these are the determinants of the sublattices of index $2^m$ of the critical lattices of the corresponding  balls $D_p$.
  
\begin{remark}
  Proposition \ref{p7} is a strengthening of Proposition \ref{p2} and its extension on  domains.  
\end{remark}

\begin{theorem}
\label{opd}
The optimal lattice packing of the Minkowski, Davis, and Chebyshev-Cohn domains is realized with respect to the sublattices of index $2^m$ of the critical lattices
  \[
  (1,0)\in\Lambda_{p}^{(0)},\; (-2^{-1/p},2^{-1/p}) \in \Lambda_{p}^{(1)}.
  \]
    \end{theorem}
{\bf  Proof}.
 By Proposition \ref{p2} the critical lattice of $2^m D_p$ is  the sublattice  of index $2^m$ of the critical lattice 
of the ball $D_p$ .
So it is the admissible lattice for $2^m D_p$ and by Proposition \ref{p1} 
is packing lattice of $2^{m-1} D_p$.
By Proposition \ref{p3} the corresponding lattice packing has maximal density and so is optimal.

\section{Inscribed and circumscribed hexagons of minimum areas}

Denote by $\Delta (2^m D_p)$ the critical determinant of the domain $2^m D_p$.
By formulas (\ref{cdde0},\ref{cdde1}) we have:
\begin{equation}
\label{cdd}
\Delta(2^m D_p) = \left\{
                   \begin{array}{lc}
     {\Delta^{(0)}_p}(2^m D_p), \; 1 \le p \le 2, \; p \ge p_{0},\\
    {\Delta^{(1)}_p}(2^m D_p), \;  2 \le p \le p_{0};\\
                     \end{array}
                       \right.
 \end{equation}
here $p_{0}$ is a real number that is defined unique by conditions
$\Delta(p_{0},\sigma_p) = \Delta(p_{0},1),  \;
2,57 < p_{0}  < 2,58, \; p_0  \approx 2.5725 $\\

Denoted by  $Ihma_{2^m D_p}$ the minimal area of  hexagons which inscibed in the domain $2^m D_p$  and have three pairs of points on the boundary of $2^m D_p$. 
From Theorems \ref{ggm}, \ref{opd} and \cite{lek} we have
\begin{theorem}
\label{ihma}
\[
  Ihma_{2^m D_p} = 3 \cdot \Delta (2^m D_p).
\]
 \end{theorem}
 
 Respectively denoted by  $Shma_{2^m D_p}$ the minimal area of  hexagons which circumscribed to the domain $2^m D_p$  and have three pairs of points on the boundary of $2^m D_p$. 
From Theorems \ref{ggm}, \ref{opd} and \cite{lek} we have

\begin{theorem}
\label{chma}
\[
        Shma_{2^m D_p}   = 4 \cdot \Delta (2^m D_p).
        \]
 \end{theorem}


\section{ Direct systems and  direct limits}
\label{secdsdld}

Various variants of direct systems and direct limits are introduced and studied in \cite{Pontryagin,bur,stei,am}.\\
Here we give their simplified versions sufficient for our purposes.
Let $X$ be a set. By a binary relation over $X$ we understand a subset of the Cartesian product $X \times X$.
By a preorder on $X$ we understand a binary relation over $X$ that is reflexive and transitive.

\begin{definition}
\label{dset}
 Let ${\mathbb N}_0$ be the set of natural numbers with zero.
A preoder $N$ on ${\mathbb N}_0$  is called a {\it directed set} if for each pair $k, m \in N$ there exists an $n \in N$ for which 
$k \le n$ and $m \le n$.
A subset  $N'$ is {\it cofinal} in  $N$ if, for each $m \in N$ their exists an $n \in N'$ such that $m \le n$.
\end{definition}

\begin{definition}
\label{dds}
A {\it direct} (or {\it inductive}) {\it system of sets} $\{X, \pi\}$ over a directed set $N$ is a function which attaches to each 
$m \in N$ a set $X^m$, and to each pair $m, n$  such that $m \le n$ in $N$, a map $\pi_m^n: X^m \to X^n$ such that, for each 
$m \in N$,  $\pi_m^m = Id$, and for $m \le n \le k$ in $N$, $\pi_n^k \pi_m^n = \pi_m^k$.
\end{definition}

\begin{remark}
We will consider direct systems of sets, topological spaces, groups and free ${\mathbb Z}$-modules.
\end{remark}

 \begin{proposition}
A directed set $N$ forms a category with elements are natural numbers $n, m, \ldots$  with zero and with morphisms
 $m \to n$ defined by the relation $m \le n$.~ A direct system over $N$ is a covariant functor from $N$ to the category of sets and maps, or to the category of topological spaces and continuous mappins, or to the category of groups  and homomorphisms or to the category of  ${\mathbb Z}$-modules and homomorphisms.
 \end{proposition}
 {\bf  Proof}. Obviously.

\begin{definition}
\label{dsm}
Let $\{X, \pi\}$ and $\{Y, \psi\}$ be direct systems over $M$ and $N$ respectively. Then a map 
\[
  \Phi:   \{X, \pi\}  \to \{Y, \psi\}
\] 

consisting of a map $\phi: M \to N$, and for each $m \in M$, a map
\[
 {\phi}^m:  X^m \to Y^{{\phi}(m)}
\]
 such that, if $m \le n$ in $M$, then commutativity holds in the diagram
\[
 \begin{CD}
   X^m @>\pi>> X^n\\
   @VV{\phi}V          @VV{\phi}V\\
   Y^{\phi(m)} @>{\psi}>>    Y^{\phi(n)}
 \end{CD}
\]

\end{definition}

\subsection{ Direct systems of Minkowski domains and their limits}

Here we will consider direct systems of Minkowski balls and domains as well as direct systems  of critical lattices.
 We use the notation according to Remark \ref{r2}.

The direct system of Minkowski balls and domains has the form (\ref{dsm}), where the multiplication by $2$  is the continuous mapping
\begin{equation}
\label{dsm}
 \begin{CD}
 D_p @>2>> 2 D_p @>2>> 2^2 D_p @>2>> \cdots @>2>> 2^m D_p @>2>> \cdots 
 \end{CD}
\end{equation}

The direct system of  critical lattices  has the form (\ref{dsml}), where the multiplication by $2$  is the homomorphism of abelian groups
\begin{equation}
\label{dsml}
 \begin{CD}
 \Lambda_{p} @>2>> 2 \Lambda_{p} @>2>> 2^2 \Lambda_{p} @>2>> \cdots @>2>> 2^m \Lambda_{p} @>2>> \cdots 
 \end{CD}
\end{equation}

Note first that real planes ${\mathbb R}^2 = (x,y)$ with Minkowski norms
\[
   |x|^p + |y|^p, \; p \in {\mathbb R}, p \ge 1,
\]
 are Banach spaces ${\mathbb B}^2_p$ .
By (general) Minkowski balls we understand balls in ${\mathbb B}^2_p$ of the form
  \begin{equation}
   \label{eq1}
   D_p: \;  |x|^p + |y|^p \le 1
  \end{equation}. 
  
 By Pontryagin duality \cite{Pontryagin} every Banach space considered in its additive structure which is an abelian topological group $G$ is reflexive that means  the existance of topological isomorphism bitween $G$ and its bidual $G^{\wedge \wedge}$.

       In our considerations we have direct systems of Minkowski balls, Minkowski domains and direct systems of critical lattices with respective maps and homomorphisms. Consider system  (\ref{dsm}) and its direct limit from the point of view of metric geometry. 
System (\ref{dsm}) is a direct system of neighborhoods of zero, and its direct limit is a 
neighborhood of zero constructed in accordance with the rules for constructing 
the direct limit of the corresponding direct system of topological spaces.
We denote this direct limit by $D^{dirlim}_p$.Accordingly, we denote the direct limit of critical lattices by $\Lambda_{p}^{dirlim}$.

\subsubsection{Direct systems}

We consider the following classes of direct systems of Minkowski balls and domains and their critical lattices.

\begin{itemize}
  
\item  {\it  Minkowski direct systems of balls and domains}:   intr $m \ge 1$, for $1 \le p<2$;
\begin{equation}
\label{dsmbd}
 \begin{CD}
 D_p @>2>> 2 D_p @>2>> 2^2 D_p @>2>> \cdots @>2>> 2^m D_p @>2>> \cdots 
 \end{CD}
\end{equation}

\item {\it  Minkowski direct systems of  critical lattices}:   int $m \ge 1$, for $1 \le p<2$;
\begin{equation} 
\label{dsmcl}
 \begin{CD}
 \Lambda_{p} @>2>> 2 \Lambda_{p} @>2>> 2^2 \Lambda_{p} @>2>> \cdots @>2>> 2^m \Lambda_{p} @>2>> \cdots 
 \end{CD}
\end{equation}

\item {\it Davis direct systems of  balls and domains}:   int $m \ge 1$, for $p_{0} > p \ge 2$;
\begin{equation}
\label{dsdbd}
 \begin{CD}
 D_p @>2>> 2 D_p @>2>> 2^2 D_p @>2>> \cdots @>2>> 2^m D_p @>2>> \cdots 
 \end{CD}
\end{equation}

\item {\it Davis direct systems of  critical lattices}:   int $m \ge 1$, for $p_{0} > p \ge 2$;
\begin{equation}
\label{dsdcl}
 \begin{CD}
 \Lambda_{p} @>2>> 2 \Lambda_{p} @>2>> 2^2 \Lambda_{p} @>2>> \cdots @>2>> 2^m \Lambda_{p} @>2>> \cdots 
 \end{CD}
\end{equation}

Below, for brevity, we will denote the balls, domains, and lattices of Chebyshev-Cohn by CC.

\item {\it CC direct systems of  balls and domains}:  int $m \ge 1$,  for $ p \ge p_{0}$;
\begin{equation}
\label{dsccbd}
 \begin{CD}
 D_p @>2>> 2 D_p @>2>> 2^2 D_p @>2>> \cdots @>2>> 2^m D_p @>2>> \cdots 
 \end{CD}
\end{equation}

\item {\it CC direct systems of  critical lattices}:  int $m \ge 1$,  for $ p \ge p_{0}$;
\begin{equation}
\label{dscccl}
 \begin{CD}
 \Lambda_{p} @>2>> 2 \Lambda_{p} @>2>> 2^2 \Lambda_{p} @>2>> \cdots @>2>> 2^m \Lambda_{p} @>2>> \cdots 
 \end{CD}
\end{equation}

\end{itemize}

\subsubsection{Direct limits}

Let us calculate direct limits of these derect systems.
Let ${\mathbb Q}_2$ and ${\mathbb Z}_2$ be respectively the field of $2$-adic numbers and its ring of integers.

 \begin{proposition}
  $D^{dirlim}_p = \varinjlim  2^m D_p  \in ({\mathbb Q}_2 / {\mathbb Z}_2)  D_p =
   (\bigcup_m \frac{1}{2^m} {\mathbb Z}_2/{\mathbb Z}_2) D_p.$
 \end{proposition}
 
  {\bf  Proof}. Follow from properties of direct systems and their direct limits \cite{bur,stei}.

 \begin{proposition}
  $\Lambda_{p}^{dirlim} = \varinjlim  2^m \Lambda_{p}  \in ({\mathbb Q}_2 / {\mathbb Z}_2)  \Lambda_{p} = 
  (\bigcup_m \frac{1}{2^m} {\mathbb Z}_2/{\mathbb Z}_2) \Lambda_{p}.$
 \end{proposition}
 
 {\bf  Proof}. Follow from properties of direct systems and their direct limits \cite{bur,stei}.\\

\end{document}